\documentclass{gtart_h}

\def\ifplaintex{\expandafter\ifx\csname documentclass\endcsname\relax}

\def\gtp{{\mathsurround=0pt\it $\cal G\mskip-2mu$eometry \&\ 
$\cal T\!\!$opology $\cal P\!$ublications}}  % GT publications

\def\recd{{\small Received:\qua\receiveddate\ifx\reviseddate\relax
\else\qquad Revised:\qua\reviseddate\fi\par}} 

\def\lognumber#1{\def\thelognumber{#1}}
\def\volumenumber#1{\def\thevolumenumber{#1}}
\def\volumeyear#1{\def\thevolumeyear{#1}}
\def\papernumber#1{\def\thepapernumber{#1}}
\def\pagenumbers#1#2{\def\startpage{#1}\def\finishpage{#2}}
\def\published#1{\def\publishdate{#1}}

\def\received#1{\def\receiveddate{#1}}

\def\accepted#1{\def\accepteddate{#1}}

\long\def\asciiabstract#1{\long\def\theasciiabstract{#1}}

\let\\\par\let\thelognumber\relax\let\thevolumenumber\relax
\let\thepapernumber\relax\let\thevolumeyear\relax\let\startpage\relax
\let\finishpage\relax\let\publishdate\relax\let\receiveddate\relax
\let\reviseddate\relax\let\accepteddate\relax\let\theasciititle\relax
\let\theasciiauthors\relax
\let\theasciiabstract\relax

\let\theasciiemail\relax

\ifplaintex
\font\logobig=cmssbx10 scaled 3836
\font\logomed=cmssbx10 scaled 2557
\else
\font\logobig=cmssbx10 scaled 4200
\font\logomed=cmssbx10 scaled 2800
\fi

\long\def\makeagttitle{   %%% start of definition of \makeagttitle
\count0=\startpage
\agt\hfill      %   Journal title (top left) 
\hbox to 45truept{\vbox to 0pt{\vglue -13truept{\logomed A\kern -.37em{\logobig 
T}\kern -.38em G}\vss}\hss}
\break
{\small Volume \thevolumenumber\ (\thevolumeyear)
\startpage--\finishpage\nl
Published: \publishdate}

\vglue .25truein

{\parskip=0pt\leftskip 0pt plus
1fil\def\\{\par\smallskip}{\Large\bf\thetitle}\par\medskip} \vglue
0.05truein

{\parskip=0pt\leftskip 0pt plus 1fil\def\\{\par}{\sc\theauthors}
\par\medskip}%
 
\vglue 0.03truein 

{\small\leftskip 25truept\rightskip 25truept{\bf Abstract}\stdspace\theabstract

{\bf AMS Classification}\stdspace\theprimaryclass
\ifx\thesecondaryclass\relax\else; \thesecondaryclass\fi\par
{\bf Keywords}\stdspace \thekeywords\par}\vglue 7truept

}   %%%% end of definition of \makeagttitle

\ifplaintex
\font\phead=cmsl9 scaled 950
\font\pnum=cmbx10 scaled 913
\font\pfoot=cmsl9 scaled 950
\headline{\vbox to 0pt{\vskip -4.5mm\line{\small\phead\ifnum
\count0=\startpage ISSN 1472-2739 (on-line) 1472-2747 (printed)
\hfill {\pnum\folio}\else\ifodd\count0\def\\{ }% 
\ifx\theshorttitle\relax\thetitle\else\theshorttitle\fi\hfill{\pnum\folio}
\else\def\\{ and }{\pnum\folio}\hfill\ifx\theshortauthors\relax\theauthors
\else\theshortauthors\fi\fi\fi}\vss}}
\footline{\vbox to 0pt{\vglue 0mm\line{\small\pfoot\ifnum\count0=\startpage
\copyright\ \gtp\hfill\else
\agt, Volume \thevolumenumber\ (\thevolumeyear)\hfill\fi}\vss}}
\else
\font=cmsl9 scaled 1050
\font\lnum=cmbx10 
\font=cmsl9 scaled 1050
\makeatletter
\def\@oddhead{{\small\lhead\ifnum\count0=\startpage ISSN 1472-2739 
(on-line) 1472-2747 (printed)\hfill {\lnum\number\count0}\else\ifodd\count0
\def\\{ }\ifx\theshorttitle\relax \thetitle \else\theshorttitle\fi\hfill
{\lnum\number\count0}\else\def\\{ and }{\lnum\number\count0}
\hfill\ifx\theshortauthors\relax 
\theauthors\else\theshortauthors\fi\fi\fi}}\def\@evenhead{\@oddhead}
\def\@oddfoot{\small\lfoot\ifnum\count0=\startpage\copyright\ \gtp\hfill\else
\agt, Volume \thevolumenumber\ (\thevolumeyear)\hfill\fi}
\def\@evenfoot{\@oddfoot}
\makeatother
\fi
\let\maketitlepage\makeagttitle
\let\makeshorttitle\maketitlepage
\let\maketitle\maketitlepage

\newwrite\gtoutfile
\long\gdef\makeheadfile{  %%% start of definition of \makeheadfile
{\def\\{, }\def\s{ }
\immediate\openout\gtoutfile head.xxx
\immediate\write\gtoutfile{Proxy-for: \ifx\theasciiauthors\relax
\theauthors\else\theasciiauthors\fi\s<\ifx\theasciiemail\relax\theemail\else\theasciiemail\fi>}
\immediate\write\gtoutfile{\noexpand\\}
\immediate\write\gtoutfile{Authors: \ifx\theasciiauthors\relax
\theauthors\else\theasciiauthors\fi}
{\def\\{ }\immediate\write\gtoutfile{Title: \ifx\theasciititle\relax
\thetitle\else\theasciititle\fi}}
\immediate\write\gtoutfile{Subj-class: GT or SG, GR etc}
\immediate\write\gtoutfile{MSC-class: \theprimaryclass\ifx\thesecondaryclass\relax\else, \thesecondaryclass\fi}
\immediate\write\gtoutfile{Journal-ref: Algebraic and Geometric Topology \thevolumenumber\s
(\thevolumeyear) \startpage-\finishpage}
\immediate\write\gtoutfile{Comments: Published by Algebraic and
Geometric Topology at}
\immediate\write\gtoutfile{\s\s\s  http://www.maths.warwick.ac.uk/agt/AGTVol\thevolumenumber/agt-\thevolumenumber-\thepapernumber.abs.html}
\immediate\write\gtoutfile{\noexpand\\}
\immediate\write\gtoutfile{}
\ifx\theasciiabstract\relax
\immediate\write\gtoutfile{\theabstract}\else
\immediate\write\gtoutfile{\theasciiabstract}\fi
\immediate\write\gtoutfile{}
\immediate\write\gtoutfile{\noexpand\\}
\immediate\write\gtoutfile{}
\immediate\closeout\gtoutfile}}  %%% end of definition of \makeheadfile

\def\maketitlepage{\makeagttitle\makeheadfile}
\let\makeshorttitle\maketitlepage
\let\maketitle\maketitlepage

\lognumber{5}
\volumenumber{4}
\volumeyear{2004}
\papernumber{5}
\published{14 February 2004}
\pagenumbers{73}{80}
\received{7 January 2004}
\accepted{19 January 2004}

\usepackage{rlepsf,amssymb, amsmath}

\newtheorem{thm}{Theorem}[section]
\newtheorem{lem}[thm]{Lemma}
\newtheorem{cor}[thm]{Corollary}

\def\dfn#1{{\em #1}}

\def\Z{\text{$\mathbb{Z}$}}

\def\C{\text{$\mathbb{C}$}}

\begin{document}
\title{On symplectic fillings}

\author{John B Etnyre}
\address{Department of Mathematics,
University of Pennsylvania\\209 South 33rd St,
Philadelphia, PA 19104-6395, USA}
\email{etnyre@math.upenn.edu}
\urladdr{http://math.upenn.edu/~etnyre}

\begin{abstract}
In this note we make several observations concerning symplectic
fillings. In particular we show that a (strongly or weakly)
semi-fillable contact structure is fillable and any filling embeds as
a symplectic domain in a closed symplectic manifold. We also relate
properties of the open book decomposition of a contact manifold to its
possible fillings. These results are also useful in showing the
contact Heegaard Floer invariant of a fillable contact structure does
not vanish \cite{OzsvathSzabo} and property P for knots
\cite{KronheimerMrowka}.
\end{abstract}

\asciiabstract{% 
In this note we make several observations concerning symplectic
fillings.  In particular we show that a (strongly or weakly)
semi-fillable contact structure is fillable and any filling embeds as
a symplectic domain in a closed symplectic manifold.  We also relate
properties of the open book decomposition of a contact manifold to its
possible fillings. These results are also useful in proving property P
for knots [P Kronheimer and T Mrowka, Geometry and Topology, 8 (2004)
295-310] and in showing the contact Heegaard Floer invariant of a
fillable contact structure does not vanish [P Ozsvath and Z Szabo,
Geometry and Topology, 8 (2004) 311-334].}

\maketitle

\section{Introduction}
It is important to consider various notion of ``symplectic filling'' when studying contact and
symplectic manifolds (especially in low dimensions). In particular, it is useful to understand when
a contact structure on a three manifold has a symplectic filling, what type of filling is it and if
there are restrictions on the topology of the filling.

Consider a contact three manifold $(M,\xi)$ and symplectic manifold $(X,\omega)$ with $M$ a boundary
component of $X$ and $\omega|_\xi>0.$ If all the boundary components of $X$ are convex then
$(X,\omega)$ is a weak (or possibly strong) semi-filling of $(M,\xi).$ (Throughout this paper all
contact structures will be (co)oriented, contact manifolds will be oriented by their contact
structures and symplectic manifolds will be oriented by their symplectic structures.) We say
the boundary component $M$ of $X$ can be \dfn{symplectically capped off} if there is a symplectic
manifold $(X',\omega')$ such that $X'\setminus M= X\cup C, \omega'|_X=\omega$ and $C$ is compact.
(Note this implies that $\partial C = -M.$) The symplectic manifold $(C,\omega'|_C)$ is called a
symplectic cap for $(M,\xi)\subset (X,\omega).$ (The notation $A\setminus B$ means the union of the
metric completions of the components of $A$ minus $B.$) Many of our results rely on the following
theorem.
\begin{thm}\label{main}
If $(X,\omega)$ is a (weak or strong) symplectic filling (or semi-filling) of $(M,\xi)$ then the
boundary component $M$ of $X$ can be capped off symplectically.
\end{thm}
We have the following immediate corollaries.
\begin{cor}\label{semi=h}
A weakly, respectively strongly, semi-fillable contact structure is weakly, respectively strongly,
fillable. \qed
\end{cor}
\begin{cor}\label{canembed}
Any symplectic four manifold with weakly or strongly convex boundary components can be embedded as a
domain in a closed symplectic manifold.\qed
\end{cor}
The theorem and both corollaries have also been established by Eliashberg in \cite{Eliashberg03}.
For {\em strong} (semi)fillings these results were essentially proven by Gay in \cite{Gay02} and
follow from Theorem~1.3 in \cite{EtnyreHonda02a}. Moreover, after writing this paper the author
learned that Stipsicz and Ghiggini also know how to prove these results for strong fillings.
Moreover, in \cite{Stipsicz}, Stipsicz has a different proof of Lemma~\ref{bord} (below), which is
the key ingredient in the proof of Theorem~\ref{main}. In \cite{LiscaMatic97} Lisca and Mati\'c and
in \cite{AO} Akbulut and Ozbagci proved the corollaries for Stein fillings. Not only do these
corollaries illuminate the nature of symplectic and contact geometry in low dimensions they have
many important applications. Specifically, Kronheimer and Mrowka \cite{KronheimerMrowka} can use
Corollary~\ref{semi=h} to prove the Property P conjecture for knots and Ozsv\'ath and Szab\'o
\cite{OzsvathSzabo} can use Corollary~\ref{canembed} to shown that their Heegaard Floer contact
invariant does not vanish for fillable contact structures. Moreover, they can use the corollary to
give an alternate proof that if $p$ surgery on a knot $K$ yields a manifold which is (orientation
preserving) homeomorphic to $L(p,1)$ then $K$ is the unknot (this result first appeared in
\cite{KMOS} using Seiberg--Witten Floer Homology).

We now turn to the topology of a symplectic filling. A fundamental result of Giroux \cite{Giroux??}
(see also \cite{Goodman}) says that a contact structure $\xi$ on a three manifold $M$ is always
supported by an open book $(L,\phi).$ That is, comes from the construction of Thurston and
Winkelnkemper \cite{TW}; or more explicitly, there is a link $L$ in $M$ that is transverse to $\xi$
such that $M\setminus L$ is fibered by Seifert surfaces for $L$ and the contact structure can be
isotoped to be arbitrarily close to the fibers in the fibration (while keeping $L$ transverse). The
\dfn{monodromy} of the open book $\phi$ is the monodromy of the fibration $M\setminus L.$
\begin{thm}\label{restrict}
Suppose the contact structure $\xi$ on a three manifold $M$ is supported by the open book $(L,\phi)$
and there is a diffeomorphism $\psi$ of the fiber such that $(L,\psi\circ\phi)$ is an open book for
$S^3$ and $\psi$ can be written as a composition of right handed Dehn twists. If $(X,\omega)$ is a,
weak or strong, symplectic semi-filling of $(M,\xi)$ then $(X,\omega)$ is a filling of $(M,\xi).$
Moreover, if $M$ is a rational homology sphere then $H_1(X;\Z)$ finite (and trivial if $M$ is an
integral homology sphere) and the intersection form on $H_2(X;\Z)$ is negative definite.
\end{thm}
Lisca (in \cite{Lisca98, Lisca99}) and Ohta--Ono (in \cite{OO}) have obtained a similar result under
the assumption that $M$ admits a metric of positive scalar curvature (but no assumption on the
supporting open book). We note that if $\phi^{-1}$ can be written as a product of right handed Dehn
twists then the conditions of the theorem are satisfied (see the proof of Lemma~\ref{bord}). If
$\phi^{-1}$ can be so written then $\phi$ can be written as product of left handed Dehn twists. One
might expect that this implies the associated contact structure is overtwisted, but this is not
necessarily true, only if the open book ``destabilizes.'' In general it is quite subtle to decide if
an open book destabilizes. It is also difficult to make general statements about representing a
diffeomorphism in terms of Dehn twists.
\begin{thm}\label{embed}
Under the hypothesis of Theorem~\ref{restrict} any symplectic filling of $(M,\xi)$ can be embedded
as a symplectic domain in a rational surface.
\end{thm}
While these theorems seem to provide potentially quite useful insights into the monodromy of open
books we content ourselves here with the following simple applications.
\begin{cor}
Let $\phi$ be any monodromy associated to any tight contact structure coming from perturbing a
Reebless foliation into a contact structure as in \cite{ET}. Then $\phi$ cannot be expressed as a
composition of left handed Dehn twists. Moreover, there is no diffeomorphism $\psi$ such that
$\phi\circ\psi$ is the monodromy for an open book of $S^3$ and $\psi$ is a composition of right
handed Dehn twists.
\end{cor}
As a specific example any monodromy map associated to any tight contact structure on $T^3$ cannot be
expressed as the composition of left handed Dehn twists.
\begin{proof}
As shown in \cite{ET}, such a contact structure on $M$ is semi-filled by a symplectic form on
$M\times[0,1].$ This contradicts the conclusion of Theorem~\ref{restrict} and thus the contact
structure cannot satisfy the hypothesis of Theorem~\ref{restrict}.
\end{proof}

{\bf Acknowledgments}\qua I thank Yasha Eliashberg for sharing an early version of \cite{Eliashberg03}
with me which made me aware of the many interesting applications of symplectic caps. I also thank
the referee who provided invaluable suggestions for improvements on the paper. In addition, I
gratefully acknowledge support from an NSF Career Grant (DMS--0239600) and FRG-0244663.

\section{Symplectic fillings}
We recall the various notions of symplectic fillings of contact manifolds. (For more details, see
the survey paper \cite{Etnyre98}, and for further discussion of the notions of fillability see
\cite{Eliashberg03}.)

A symplectic manifold $(X,\omega)$ is said to have \dfn{strongly convex boundary} if there is a
vector field $v$ defined in the neighborhood of $\partial X$ that points transversely out of $X$ and
dilates $\omega$ ({\em ie}, $\mathcal{L}_v\, \omega=\omega$). The form
$\alpha=(\iota_v\, \omega)\vert_{\partial X}$ is a contact form on $\partial X.$ The contact manifold
$(M=\partial X,\xi)$ on the boundary is said to be \dfn{strongly symplectically fillable} by
$(X,\omega).$ If $(X,\omega)$ supports a Stein structure then $(M,\xi)$ is said to be
\dfn{holomorphically fillable} (this is sometimes called Stein fillable). A symplectic manifold $(X,
\omega)$ is said to have \dfn{weakly convex boundary} if $\partial X$ admits a contact structure
$\xi$ such that $\omega\vert_\xi >0.$ A contact manifold $(M,\xi)$ is \dfn{weakly symplectically
fillable} if it is the weakly convex boundary of a symplectic manifold. A contact structure is said
to be weakly (or strongly) \dfn{semi-fillable} if it is one component of the boundary of a weakly
(or strongly) convex symplectic manifold. (Note a Stein manifold always has a connected boundary so
there is no notion of holomorphically semi-fillable.) When talking about strong or weak symplectic
fillings we will frequently drop the word ``symplectic'' and just refer to strong or weak fillings.

The following diagram indicates the hierarchy of contact structures.
$$
\begin{array}{ccc}
& &\\
\mbox{Tight} & &   \\
\cup \not | & &\\
\mbox{Weakly semi-fillable} & \supsetneqq & \mbox{Strongly
semi-fillable}\\ ||  & & || \\
\mbox{Weakly fillable} & \supsetneqq  & \mbox{Strongly
fillable}\\ & & \cup\\
& & \mbox{Holomorphically fillable}\\
& &\\
\end{array}
$$
The properness of the inclusion of the set of weakly symplectically semi-fillable contact structures
into the set of tight contact structures was shown in \cite{EtnyreHonda02b} (for further examples
see \cite{LiscaStipsicz}). The properness of the inclusion of the set of strongly fillable
structures into the set of weakly fillable structures is due to Eliashberg \cite{Eliashberg96} (for
further examples see \cite{DG}). The equalities of the semi-fillings and ``honest'' fillings is the
content of Corollary~\ref{semi=h} and \cite{Eliashberg03}. Though there are strong fillings of a
contact structure that are not holomorphic fillings (because there is more than one boundary
component \cite{McD}) it is not currently known if strongly fillable contact structures are also
holomorphically fillable.

If $M$ is a rational homology sphere, a weak symplectic filling can be modified into a strong
symplectic filling \cite{OO} (the germ of the argument originally appeared in \cite{Eliashberg91},
see also \cite{Eliashberg03}). If $M$ has a positive scalar curvature metric, then moreover a
semi-filling is automatically a one boundary component filling, and all four notions of symplectic
filling become the same (this is due to Lisca \cite{Lisca98} and Ohta--Ono \cite{OO}). The proof
builds on the work of Kronheimer and Mrowka \cite{KM} and relies on Seiberg--Witten theory.

\section{Proof of Theorem~\ref{main}}
For the proof of Theorem~\ref{main} we will need the following lemma (see also \cite{Stipsicz}).
\begin{lem}\label{bord}
Under the hypothesis of Theorem~\ref{main}, 2--handles may be attached to $X$ to form a new
symplectic manifold with weakly convex boundary $(X',\omega')$ such that $X'\setminus M = X \cup B$
with $\partial B= M'\cup (-M)$ where $M'$ is a homology sphere.
\end{lem}
With this lemma in hand we prove Theorem~\ref{main}
\begin{proof}[Proof of Theorem~\ref{main}]
Given $(M,\xi)$ and $(X,\omega)$ as in the theorem, apply the above lemma to obtain $(X',\omega').$
Since $M'$ is a homology sphere the results of Ohta--Ono mentioned in the previous section imply
$\omega'$ may be isotoped near $M'$ so that $(X',\omega')$ is strongly convex along $M'.$
Theorem~1.3 from \cite{EtnyreHonda02a} says that any contact manifold has (infinitely many) strongly
concave fillings. Since a strongly convex and strongly concave filling may be glued to form a
symplectic manifold \cite{Etnyre98} we see that $(X',\omega')$ may be symplectically capped off. Let
$C$ be the cap. Then $C\cup B$ ($B$ from Lemma~\ref{bord}) is the symplectic cap for $(M,\xi).$
\end{proof}

We are left to prove the lemma.
\begin{proof}[Proof of Lemma~\ref{bord}]
Giroux's result about open books and contact structures implies we can assume that $(M,\xi)$ is
supported by an open book with connected binding (this can be achieved by ``positive Hopf
stabilization''). Let $\phi$ be the monodromy of the open book. Since the mapping class group is
generated by right handed Dehn twists about arbitrary loops plus arbitrary Dehn twists parallel to
the boundary, we can write $\phi^{-1}$ as $\phi_0\circ \phi_1,$ where $\phi_0$ is a composition of
right handed Dehn twists and $\phi_1$ is some number of left handed Dehn twists parallel to the
boundary (we can put all the boundary parallel twists at the end since they commute with the other
Dehn twists). Thus by attaching 2--handles to $X$ along $M$ in a symplectic way (see
\cite{Eliashberg90, Gompf98, Weinstein91}) we can get a symplectic manifold with a boundary
component $(M'',\xi''),$ such that $(M'',\xi'')$ is supported by an open book with monodromy
$\phi\circ \phi_0$ which is, say, $k$ right handed Dehn twists parallel to the boundary.
Topologically the manifold is shown in the Figure~\ref{firstm}.
\begin{figure}[ht!]
  \relabelbox \small {\epsfxsize=4.5in\centerline{\epsfbox{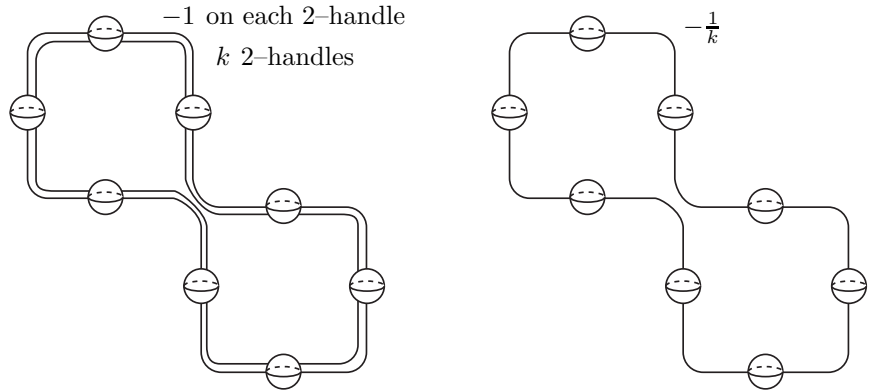}}} 
  \relabel{-1 on each 2-handle}{$-1$ on each 2--handle}
  \relabel {k 2-handles}{$k$ 2--handles}
  \relabel {-1/k}{$-\frac1k$}
  \endrelabelbox
        \caption{On the left hand side we show an open book with $k$ right handed Dehn twists
          parallel to the boundary. The right hand side is homeomorphic to the left hand side.}
        \label{firstm}
\end{figure}
\iffalse\begin{figure}[ht!]\small  
  \psfrag {-1 on each 2-handle}{$-1$ on each 2--handle}
  \psfrag {k 2-handles}{$k$ 2--handles}
  \psfrag {-1/k}{$-\frac1k$}
  {\centerline{\includegraphics[width=4.5in]{openbook.eps}}}
        \caption{On the left hand side we show an open book with $k$ right handed Dehn twists
          parallel to the boundary. The right hand side is homeomorphic to the left hand side.}
        \label{firstm}
\end{figure}\fi
Now for each 1--handle in the figure add a 2--handle that runs over it once (this can be done in a
symplectic way). We now have a symplectic manifold with boundary $(M',\xi').$ Topologically $M'$ is
obtained by $-\frac1k$ surgery on a knot in the three sphere, thus it is a homology sphere.
\end{proof}

\section{Topology of symplectic fillings}
\begin{proof}[Proof of Theorem~\ref{restrict} and \ref{embed}]
Given the hypothesized diffeomorphism $\psi$ we can repeat the proof of Lemma~\ref{bord} and attach
2--handles to $X$ along $M$ to form a new symplectic manifold $(X',\omega')$ with new boundary
component $(S^3,\xi').$ Since $\xi'$ is fillable we know it is tight. Thus $\xi'$ is the unique
tight contact structure on $S^3$ (see \cite{Eliashberg92}). McDuff's well known extension
\cite{McDuff} of Gromov's theorem \cite{Gromov85} says that $(X',\omega')$ is symplectomorphic to
the blow up of the standard symplectic structure on $B^4\subset \C^2.$ This of course implies our
original symplectic filling $(X,\omega)$ only had one boundary component and that we may cap off
$(X',\omega')$ by a symplectic disk bundle over $S^2$ to obtain a rational surface. This completes
the proof of Theorem~\ref{embed}. To finish the proof of Theorem~\ref{restrict} simply note that if
a rational homology sphere splits a rational surface as ours does then the component $(X',\omega')$
must be negative definite and $H_1(X;\Z)$ must be torsional (and trivial if $M$ is an integral
homology sphere).
\end{proof}

\Addresses\recd
\end{document}